\documentclass[12pt]{amsart}
\pdfoutput=1 
\usepackage{hyperref}

\usepackage{amsmath,amsfonts,amssymb}
\usepackage{mathrsfs}
\usepackage{verbatim}
\usepackage{amsthm}

\usepackage[english]{babel}

\newtheorem{theorem}{Theorem}[section]
\newtheorem{lemma}[theorem]{Lemma}

\newtheorem{corollary}[theorem]{Corollary}

\def \bH {\mathbb H}

\def \bN {\mathbb N}

\def \bR {\mathbb R}

\def \cF {\mathcal F}

\def \cL {\mathcal L}

\def \cR {\mathcal R}
\def \cS {\mathcal S}

\def \fh {\mathfrak h}

\def \tr {\text{\rm Tr}}

\def \ad {\text{\rm ad}}
\def \id {\text{\rm I}}

\def\HS{{\mathtt{HS}}}

\def \Op {\text{\rm Op}}
\def \Opw {\text{\rm Op}^W}

\def \sgn {\text{\rm sgn}}

\def \FR {\cF_{\bR^{2n+1}}}

\begin{document} 

\title{A pseudo-differential calculus\\ on the Heisenberg group}

\author{V\'eronique Fischer and Michael Ruzhansky}
\email{v.fischer@imperial.ac.uk}
\email{m.ruzhansky@imperial.ac.uk}
\address{Department of Mathematics,
Imperial College London\\
180 Queen's Gate, \\
London SW7 2AZ, United Kingdom}
\maketitle

\begin{abstract}
In this note we present a symbolic pseudo-differential calculus 
on the Heisenberg group.
We particularise to this group
our general construction
 \cite{RF-monograph,RF-aalto,RF-cras1}
 of pseudo-differential calculi
on graded groups.
The relation between the Weyl quantization and the representations 
of the Heisenberg group
enables us to consider here scalar-valued symbols.
We find that the conditions defining the symbol classes 
 are similar but different to the ones in  \cite{BFG}.
 Applications are given to Schwartz hypoellipticity and to subelliptic estimates on
the Heisenberg group.


{\it MSC classes: 35S05, 43A80}
\end{abstract}
\section{Introduction}
\label{Intro}
In \cite{RF-monograph}, see also \cite{RF-aalto,RF-cras1}, 
a pseudo-differential calculus is developed
in the setting of graded Lie groups using their representations.
Here we present the results of this construction in the particular case 
of the Heisenberg group $\bH_{n}$.

It is well known that
the representations of $\bH_n$
are intimately linked with the Weyl quantization on $\bR^n$
(see e.g., \cite{Tnma}, and Section \ref{sec_schrodinger_rep+Wquantization} below).
Together with the analogue of the Kohn-Nirenberg quantization on Lie groups 
(see e.g., \cite{Tnma,RTb,RF-monograph}, and 
Section \ref{sec_quantization+symbol_classes} below),
this link enables 
the development 
of pseudo-differential calculi
on $\bH_n$
with scalar-valued symbols that depend on parameters.
However,  the remaining difficulty lies in finding 
conditions to be imposed on those symbols 
so that the resulting class of operators has the expected properties of a calculus.

Although M. Taylor explained these general ideas in the
setting of the Heisenberg groups in \cite{Tnma},
he chose to restrict his analysis in \cite{Tnma} mainly 
to invariant (i.e. convolution) operators on $\bH_n$
with symbols defined by some asymptotic expansions.
To the authors' knowledge, 
the only study of non-invariant calculi 
with scalar-valued symbols
on $\bH_n$  was done, until now,
 by H. Bahouri, C. Fermanian-Kammerer and I. Gallagher
 in \cite{BFG}.
 Their work is devoted  to the case of $\bH_n$ only.
 Moreover,  the conditions 
imposed on the scalar-valued symbols
might appear difficult to apprehend  for some readers,
as they come from  technical parts of the proofs of the calculi's properties
(see the more recent version of \cite{BFG} on the server Hal).
Our conditions on symbol classes differ from those in \cite{BFG} for small
$\lambda$.  At the end of this note we list several applications of the analysis
in our classes, to the hypoellipticity properties and subelliptic
estimates for several operators on the
Heisenberg group.
 
Our approach to find the conditions on the symbols 
is different from \cite{Tnma} and \cite{BFG}:
we particularise to the setting of $\bH_n$
our definition of pseudo-differential calculi
valid on a large class of nilpotent Lie groups, 
namely the graded groups, 
see \cite{RF-monograph,RF-aalto,RF-cras1}.
In our general construction, 
the symbols are operator-valued.
Nonetheless on $\bH_n$, 
using the link between the Weyl quantization and the representations of $\bH_n$, 
this is equivalent to using the scalar-valued symbols.
The purpose of this note is to present 
what the general conditions on the symbols given in 
\cite{RF-monograph,RF-aalto,RF-cras1}
become when expressed on the level of scalar-valued symbols of $\bH_n$.
In particular, 
we find  conditions which are similar but different to the ones in  \cite{BFG}.
As applications for our analysis, 
we give sufficient condition for Schwartz hypoellipticity and 
for subelliptic estimates on
the Heisenberg group.

\section{Schr\"odinger representations and Weyl quantization}
\label{sec_schrodinger_rep+Wquantization}

We start by fixing the notation required for presenting our results.
We realise the Heisenberg group $\bH_{n}$
as the manifold $\bR^{2n+1}$ endowed with the  law
$$
(x,y,t)(x',y',t')=(x+x',y+y',t+t'+\frac 12 (xy'-x'y)),
$$
where $(x,y,t)$ and $(x'y',t')$ are in $\bR^{n}\times \bR^{n}\times \bR=\bR^{2n+1}$.
Here we adopt the convention that
if $x$ and $y$ are two vectors in $\bR^n$ for some $n\in \bN$,
then $xy$ denotes their standard scalar product
$$
xy=\sum_{j=1}^n x_j y_j \quad\mbox{if}\quad
x=(x_1,\ldots,x_n),\ 
y=(y_1,\ldots,y_n).
$$
The canonical basis for the Lie algebra $\fh_{n}$ of $\bH_n$ is given 
by the left-invariant 
vector fields
$$
X_j=\partial_{x_j} -\frac {y_j} 2 \partial_t,\quad 
Y_j=\partial_{y_j} +\frac {x_j} 2 \partial_t,\quad 
j=1,\ldots,n,\quad \mbox{and}\quad T=\partial_t.
$$
The canonical commutation relations are 
$$
[X_j,Y_j]=T,\quad j=1,\ldots, n,
$$
and $T$ is the centre of $\fh_{n}$. 
The Heisenberg Lie algebra is stratified via 
$\fh_{n}=V_1\oplus V_2$
where $V_1$ is linearly spanned by the $X_j$'s and $Y_j$'s,
while $V_2=\bR T$.
Therefore,
the group $\bH_n$ is naturally equipped with the family of dilations $D_r$
 given by
$$
D_r (x,y,t)=r(x,y,t)=(rx, ry,r^2t),
\quad (x,y,t)\in \bH_{n}, \ r>0.
$$
The `canonical' positive Rockland operator in this setting is 
$\cR=-\cL$, where $\cL$ is the sub-Laplacian
$$
\cL:=
\sum_{j=1}^{n} (X_j^2 +Y_j^2)
=\sum_{j=1}^{n} \left(
\left(\partial_{x_j} -\frac {y_j} 2 \partial_t\right)^2
+
\left(\partial_{y_j} +\frac {x_j} 2 \partial_t\right)^2
\right).
$$

The Schr\"odinger representations of the Heisenberg group $\bH_{n}$ 
are the  infinite dimensional unitary representations of $\bH_{n}$
(we allow ourselves to  identify unitary representations 
with their unitary equivalence classes).
Parametrised by
$\lambda\in \bR\backslash\{0\}$,
they act on $L^2(\bR^n)$.
We denote them  by $\pi_\lambda$ and realise them as
\begin{equation}
\label{eq_pilambda}
\pi_\lambda (x,y,t)h(u)
=
e^{i \lambda (t+\frac 12 xy)} 
e^{i \sqrt{\lambda}    yu} h(u+\sqrt{|\lambda|}  x),
\end{equation}
for $h\in L^2(\bR^{n})$, $u\in \bR^n$, and $(x,y,t)\in \bH_{n}$,
where we use the  convention
\begin{equation}
\label{eq_convention_sqrtlambda}
\sqrt{\lambda}:=
\sgn(\lambda) \sqrt{|\lambda|}
=
\left\{\begin{array}{ll}
\sqrt{\lambda}&\mbox{if}\ \lambda>0,\\
-\sqrt{|\lambda|}&\mbox{if}\ \lambda<0.\\
\end{array}\right.
\end{equation}
The group Fourier transform of a function $\kappa\in L^1(\bH_n)$
is by definition
$$
\widehat \kappa(\pi_\lambda)
\equiv
\pi_\lambda(\kappa)
:=
\int_{\bH_{n}} \kappa(x,y,t)
\pi_\lambda(x,y,t)^* dxdydt.
$$
As already noted in \cite{Tnma},
it can be effectively computed by
$$
\pi_\lambda(\kappa)h(u)
=
\int_{\bR^{2 n +1}} \kappa(x,y,t) 
e^{i \lambda (-t+\frac 12 xy)} 
e^{-i \sqrt{\lambda} yu} h(u- \sqrt{|\lambda|} x )
 dxdydt,
$$
for $h\in L^2(\bR^{n})$ and $u\in \bR^n$,
that is, 
\begin{equation}
\label{eq_pilamdakappa_OpW}
\widehat \kappa(\pi_\lambda)(u)
=
\pi_\lambda(\kappa)h(u)
=
(2\pi)^{\frac{n} 2}\Opw \left[\cF_{\bR^{2n+1}} (\kappa) (\sqrt {|\lambda|} \, \cdot,\sqrt \lambda\, \cdot ,\lambda)\right].
\end{equation}
Here the Fourier transform $\cF_{\bR^{2n+1}}=\cF_{\bR^N}$ is defined via
 $$
 \cF_{\bR^N} f (\xi)=
 (2\pi)^{-\frac N2} \int_{\bR^N} f(x) e^{-ix\xi} dx, 
 \quad(\xi\in \bR^N, \ f\in L^1(\bR^N)),
 $$
and $\Opw$ denotes the Weyl quantization, 
which is given for a reasonable symbol $a$
on $\bR^n\times\bR^n$, by
$$ \Opw(a) f(u) = 
  a(D,X)f(u) = 
 (2\pi)^{-n} 
 \int_{\bR^n}
 \int_{\bR^n}
 e^{i(u-v)\xi } a(\xi,\frac{u+v}2) f(v) dv d\xi,
$$
 where $f\in \cS(\bR^n)$ and $u\in \bR^n$.
We keep the same notation $\pi_\lambda$ for the 
corresponding infinitesimal representations.
We readily compute that
$$
\begin{array}{rclcl}
\pi_\lambda(X_j) 
&=& \sqrt{|\lambda|} \partial_{u_j}
&=&\Opw \left(i\sqrt{|\lambda|} \xi_j\right),
\\
\pi_\lambda(Y_j) &=&
i \sqrt{\lambda}   u_j
&=&\Opw\left( i \sqrt{\lambda}   u_j \right),
\\
\pi_\lambda(T) 
&=&i \lambda \id
&=& \Opw(i\lambda),
\end{array}
$$
thus
$$
\pi_\lambda(\cL) 
=
\sum_{j=1}^{n} (\pi_\lambda(X_j)^2 +\pi_\lambda(Y_j)^2)
=|\lambda| \sum_{j=1}^n \left(\partial_{u_j}^2 -u_j^2\right)
= -\Opw \left(|\lambda| \sum_{j=1}^n\left( \xi_j^2 + u_j^2\right)\right)
.
$$

With our choice of notation and definitions, 
the Plancherel measure is $c_n|\lambda|^{n} d\lambda$
in the  sense that the  Plancherel formula 
\begin{equation}
\label{eq_Plancherel_formula}
\int_{\bH_{n_o}} |\kappa(x,y,t)|^2 dx dy dt
=
c_n\int_{\bR\backslash\{0\}}
\|\pi_\lambda (\kappa)\|_{\HS}^2 
|\lambda|^{n} d\lambda.
\end{equation}
holds for any  $\kappa\in \cS(\bH_{n})$.
For the value of the constant $c_n$, see \cite{RF-monograph}.
Here $\|\cdot\|_{\HS}$ denotes the Hilbert-Schmidt norm of an operator on $L^2(\bR^n)$, 
that is, 
$\|B\|_{\HS}^2=\tr (B^* B)$.
This allows one to extends unitarily the definition of the group Fourier transform to $L^2(\bH_n)$. Formula \eqref{eq_Plancherel_formula} 
then holds for any $\kappa\in L^2(\bH_n)$.

\section{Difference operators}

Difference operators were defined 
 in   \cite{RTb,RTi} as acting on Fourier coefficients on compact Lie groups,  
and on graded nilpotent Lie groups in \cite{RF-monograph}.
In the setting of the Heisenberg group, this yields the definition of 
the difference  operators 
$\Delta_{x_j}$, $\Delta_{y_j}$, and $\Delta_{t}$
via
$$
\Delta_{x_j} \widehat \kappa (\pi_\lambda)
:=
\pi_\lambda(x_j \kappa),
\quad
\Delta_{y_j} \widehat \kappa (\pi_\lambda)
:=
\pi_\lambda(y_j \kappa),
\quad
\Delta_{t} \widehat \kappa (\pi_\lambda)
:=
\pi_\lambda(t \kappa),
$$
for suitable distributions $\kappa$ defined on $\bH_n$.
We can compute that
$$
\begin{array}{rclcl}
\Delta_{x_j} |_{\pi_\lambda}
&=&
\frac1{i \lambda} \ad \left(\pi_\lambda(Y_j)\right)
&=& 
\frac 1{\sqrt{|\lambda|}} \ad u_j,
\\
\Delta_{y_j} |_{\pi_\lambda}
&=&
-\frac1{i \lambda} \ad \left(\pi_\lambda(X_j)\right)
&=& -\frac 1{i\sqrt{\lambda}} \ad \partial_{u_j},
\end{array}
$$
and 
$$
\Delta_{t} |_{\pi_\lambda}
=
i \partial_{\lambda} 
+\frac 12 \sum_{j=1}^{n}\Delta_{x_j}\Delta_{y_j} |_{\pi_\lambda}
-\frac i{2\lambda}
\sum_{j=1}^{n}
\left(\pi_\lambda(Y_j)
\Delta_{y_j}|_{\pi_\lambda}   
+\Delta_{x_j} |_{\pi_\lambda} \pi_\lambda(X_j)\right).
$$
When  $\pi_\lambda(\kappa) = \Opw (a_\lambda)$ and $a_\lambda=\{a_\lambda(\xi,u)\}$,
we have
\begin{equation}
\left.\label{eq_Deltaxyt_opw}
\begin{array}{rcl}
\Delta_{x_j}\pi_\lambda(\kappa) 
&=&
\Opw \left(\frac {i}{\sqrt{|\lambda|}}   \partial_{\xi_j} a_\lambda\right)
\\
\Delta_{y_j} \pi_\lambda(\kappa) 
&=&
\Opw \left(\frac{i}{\sqrt{\lambda}}  
\partial_{u_j} a_\lambda\right)\\
\Delta_t \pi_\lambda(\kappa) &=& 
i\Opw\left(\tilde \partial_{\lambda,\xi,u}a_\lambda(\xi,u)
\right)
\end{array}\right\}
\end{equation}
where
\begin{equation}
\label{eq_def_dlambdaxiu}
\tilde \partial_{\lambda,\xi,u}:=\partial_{\lambda}
 - \frac 1{2\lambda}
 \sum_{j=1}^{n}
 \left(u_j \partial_{u_j} + \xi_j\partial_{\xi_j}\right).
\end{equation}
For example, we have
$$
\begin{array}{l}
\Delta_{x_j} \pi_\lambda(Y_k)
=\Delta_{x_j} \pi_\lambda(T)
=\Delta_{y_j} \pi_\lambda(X_k)
=\Delta_{y_j} \pi_\lambda(T)
=0,
\\
\Delta_{x_j} \pi_\lambda(X_k)
=   \Delta_{y_j} \pi_\lambda(Y_k)
=-\delta_{jk} \id,
\quad
\Delta_{t} \pi_\lambda(T)
=-\id,
 \\
\Delta_{x_j} \pi_\lambda(\cL)
= -2 \pi_\lambda(X_j),\quad
\Delta_{y_j} \pi_\lambda(\cL)
= -2\pi_\lambda(Y_j),\quad
\Delta_{t} \pi_\lambda(\cL)=0.
\end{array}
$$

\medskip

The following equalities shed some light on why, 
for example in \cite{BFG},
another normalisation of the Weyl symbol is preferred.
Indeed, the expressions on the right-hand sides in \eqref{eq_Deltaxyt_opw},
in particular for the operator $\tilde\partial_{\lambda,\xi,u}$ defined in \eqref{eq_def_dlambdaxiu},
become  very simple.

\medskip

\begin{lemma}
\label{lem}
Let $a_\lambda=\{a_\lambda(\xi,u)\}$ be a family of Weyl symbols
depending smoothly on $\lambda\not=0$.
If $\tilde a_\lambda$ is the renormalisation obtained via
$$
a_{\lambda} (\xi,u)
:=
\tilde a_{\lambda} (\sqrt{|\lambda|}\xi,\sqrt{\lambda} u),
$$
then
$$
\tilde\partial_{\lambda,\xi,u}a_\lambda(\xi,u)
=
\{\partial_\lambda 
\tilde a_{\lambda}\} (\sqrt{|\lambda|}\xi,\sqrt{\lambda} u),
$$
$$
\frac 1{i\sqrt{|\lambda|}}   \partial_{\xi_j} a_\lambda
= (\partial_{\xi_j}\tilde a_{\lambda}) (\sqrt{|\lambda|}\xi,\sqrt{\lambda} u),
\quad\mbox{and}\quad
\frac 1{i\sqrt{\lambda}}  
\partial_{u_j} a_\lambda
=
(\partial_{u_j}\tilde a_{\lambda}) (\sqrt{|\lambda|}\xi,\sqrt{\lambda} u).
$$
Consequently,
$$
\Delta_{x_j}\pi_\lambda(\kappa) 
=
i\Opw \left( \partial_{\xi_{j}}
\tilde a_{\lambda} \right),\quad
\Delta_{y_j} \pi_\lambda(\kappa) 
=
i\Opw \left(\partial_{u_j} \tilde a_\lambda\right)
\;\textrm{ and }\;
\Delta_t \pi_\lambda(\kappa) =
i\Opw\left(
\partial_{\lambda} \tilde a_\lambda
\right).
$$
\end{lemma}

\section{Quantization and symbol classes}
\label{sec_quantization+symbol_classes}

In this note, for simplicity, 
we change slightly the notation with respect to the general case developed in \cite{RF-monograph}.
Firstly we want to keep the letter $x$ to denote part of the coordinates of the Heisenberg group
and we choose to denote a general element of the Heisenberg group  by, e.g., 
$$
g=(x,y,t) \in \bH_{n}.
$$
Secondly we may define a symbol as parametrised by 
$$
\sigma(g, \lambda):=\sigma(g,\pi_\lambda),
\quad (g,\lambda)\in \bH_{n}\times \bR\backslash\{0\}.
$$
Thirdly we modify the indices $\alpha\in \bN_0^{2n+1}$
in order to write them as
$$
\alpha =(\alpha_1,\alpha_2,\alpha_3),
$$
with 
$$
\alpha_1=(\alpha_{1,1},\ldots,\alpha_{1,n}) \in \bN_0^{n},
\quad
\alpha_2=(\alpha_{2,1},\ldots,\alpha_{2,n}) \in \bN_0^{n},
\quad 
\alpha_3 \in \bN_0.
$$
The homogeneous degree of $\alpha$ is then
$$
[\alpha]=|\alpha_1|+|\alpha_2|+2\alpha_3.
$$
For each $\alpha$
we write 
$$
g^\alpha=
x^{\alpha_1}y^{\alpha_2} t^{\alpha_3} ,
\quad \mbox{where}\quad
x^{\alpha_1} = 
x_1^{\alpha_{11}}\ldots x_{n}^{\alpha_{1n}},
\quad
y^{\alpha_2} = 
y_1^{\alpha_{21}}\ldots y_{n}^{\alpha_{2n}},
$$
and we define the corresponding difference operator
$$
{\Delta'}^{\alpha} :=\Delta_{x}^{\alpha_1}\Delta_y^{\alpha_2} \Delta_t^{\alpha_3},
\quad \mbox{where}\quad
\Delta_x^{\alpha_1} := 
\Delta_{x_1}^{\alpha_{11}}\ldots \Delta_{x_{n}}^{\alpha_{1n}},
\quad
\Delta_y^{\alpha_2} := 
\Delta_{y_1}^{\alpha_{21}}\ldots \Delta_{y_{n}}^{\alpha_{2n}}.
$$
We also write $X^\alpha= X^{\alpha_1} Y^{\alpha_2}T^{\alpha_3}$, 
where 
$X^{\alpha_1} = 
X_1^{\alpha_{11}}\ldots X_{n}^{\alpha_{1n}}$,
and
$Y^{\alpha_2} = 
Y_1^{\alpha_{21}}\ldots Y_{n}^{\alpha_{2n}}$.

Following \cite{RF-monograph}, 
we define the symbol class  $S^m_{\rho,\delta}(\bH_n)$
as the set of symbols $\sigma$ 
for which all the following quantities are finite:
$$
\|\sigma\|_{S^m_{\rho,\delta},a,b,c}
:=
\sup_{\lambda\in \bR\backslash\{0\}, \, g\in \bH_{n}}
\|\sigma (g , \lambda) \|_{S^m_{\rho,\delta},a,b,c}
,\quad a,b,c\in \bN_0,
$$
where 
$$
\|\sigma (g, \lambda) \|_{S^m_{\rho,\delta},a,b,c}
:=
\!\!\!\!\!\!
\sup_{\substack{ [\alpha]\leq a\\
 [\beta]\leq b,\, |\gamma|\leq c}}
 \!\!\!\!\!\!
\|\pi_\lambda(\id-\cL)^{\frac{\rho [\alpha]-m -\delta[\beta] +\gamma}2 }
X_g^\beta{\Delta'}^\alpha \sigma(g,\lambda) 
\pi_\lambda(\id-\cL)^{-\frac{\gamma}2 }\|_{op} 
.
$$

A natural quantization 
on any type I Lie group
is the analogue of the Kohn-Nirenberg quantization on $\bR^n$,
 see, e.g., \cite{Tnma} for general remarks, \cite{RTb} for the 
 consistent development in the case of compact Lie group, and \cite{RF-monograph} for the case of nilpotent Lie groups.
In the particular case of the Heisenberg group,
this quantization 
 associates to a symbol $\sigma$ (for example in $S^m_{\rho,\delta}(\bH_n)$)
the operator $A=\Op(\sigma)$ acting on $\cS(\bH_n)$ given by
\begin{equation}
\label{eq_quantization}
A\varphi(g)= 
c_n \int_{\lambda\in \bR\backslash\{0\}} 
\tr \left(\pi_\lambda (g)\sigma(g,\lambda)\pi_\lambda(\varphi)\right) 
|\lambda|^n d\lambda .
\end{equation}
Here we have used our notation, 
especially for the Plancherel measure 
$c_n|\lambda|^n d\lambda$, 
see \eqref{eq_Plancherel_formula}. 
We denote by 
\begin{equation}
\label{eq_def_Psim}
\Psi^m_{\rho,\delta}(\bH_n)
:=\{\Op(\sigma), \sigma\in S^m_{\rho,\delta}(\bH_n)\},
\end{equation}
the class of operators corresponding to the symbols in $S^m_{\rho,\delta}(\bH_n)$ via this quantization.

\medskip

\medskip

The main result of this note shows that the symbols $\sigma=\{\sigma(g,\lambda)\}$ in $S^m_{\rho,\delta}$ are all of the form 
$\sigma(g,\lambda) 
= \Opw (a_{\lambda,g} (\xi,u))$
with $a_{\lambda,g}$ (called {\em $\lambda$-symbols}) satisfying  properties 
of Shubin type:

\medskip

\begin{theorem}
\label{thm_characterisation_Sm}
Let $\rho,\delta$ be real numbers such that 
$1\geq \rho\geq \delta\geq 0$
and $(\rho,\delta)\not=0$.
\begin{enumerate}
\item 
If $\sigma=\{\sigma(g,\lambda)\}$ is in $S^m_{\rho,\delta}(\bH_n)$ then
there exists a unique smooth function
$a=\{ a(g,\lambda,\xi,u)= a_{g,\lambda}(\xi,u)\}$ on $\bH_{n}\times \bR\backslash \{0\}\times \bR^{n}\times\bR^{n}$  
such that 
\begin{equation}
\label{eq_sigma_opw_tildea}
\sigma(g,\lambda)=\Opw\left( a_{g,\lambda} \right).
\end{equation}
It satisfies for any $\alpha,\beta\in \bN_0^n$, 
$\tilde \beta\in \bN_0^{2n+1}$ and $\tilde \alpha\in \bN_0$,
\begin{equation}
\label{eq_cond_aglambda}
|\partial_\xi^\alpha\ \partial_u^\beta \
\tilde\partial_{\lambda,\xi,u}^{\tilde \alpha} \
X^{\tilde \beta}_g \ a_{g,\lambda}(\xi,u)|
\leq 
C_{\alpha,\beta,\tilde\alpha,\tilde \beta}
|\lambda|^{\rho\frac{|\alpha|+|\beta|}2}
\left(1+|\lambda|(1+|\xi|^2+|u|^2)\right)^{\frac{m-2\rho\tilde \alpha+\delta [\tilde \beta]-\rho(|\alpha|+|\beta|)} 2},
\end{equation}
where  the operator $\tilde\partial_{\lambda,\xi,u}$ was defined in \eqref{eq_def_dlambdaxiu}.

\item 
Conversely, if 
$a=\{a(g,\lambda,\xi,u)=a_{g,\lambda}(\xi,u)\}$ is  a smooth function on 
$\bH_{n}\times \bR\backslash \{0\}\times \bR^{n}\times\bR^{n}$ 
satisfying \eqref{eq_cond_aglambda}
for every $\alpha,\beta\in \bN_0^n$, $\tilde \alpha\in \bN_0$,
then there exists a unique symbol $\sigma\in S^m_{\rho,\delta}(\bH_n)$
such that \eqref{eq_sigma_opw_tildea} holds.

\item The resulting class of operators 
$\cup_{m\in \bR} 
 \Psi^m_{\rho,\delta}(\bH_n)$
is an algebra of operators, 
the product being the usual composition.
It is stable under taking the adjoint and 
 contains the left-invariant differential calculus.
Each operator $A\in \Psi^m_{\rho,\delta}(\bH_n)$ maps continuously the Sobolev space $L^2_s(\bH_n)$ of the Heisenberg group to $L^2_{s-m}(\bH_n)$
with the loss of $m$ derivatives (for any $s\in \bR$).
\end{enumerate}
\end{theorem}

For the definition of the Sobolev spaces $L^2_s(\bH_n)$ on $\bH_n$
and more generally on any stratified group, see \cite{F-ark.75}.
Part $(iii)$ summarises the main results of 
the general construction made in \cite{RF-monograph} 
on any graded groups.

\medskip

Let us notice that 
writing $\sigma(g,\lambda)=\Opw (a_{g,\lambda})$ as in \eqref{eq_sigma_opw_tildea}
and using \eqref{eq_pilamdakappa_OpW} for 
$$
\pi_\lambda(\varphi) \pi_\lambda(g) = \pi_\lambda(\varphi(g\, \cdot))
$$ 
yield the following 
alternative formula for the quantization 
given in \eqref{eq_quantization}:
\begin{equation}\label{EQ:quant2}
A\varphi(g)= 
c'_n
\int_{\lambda\in \bR\backslash\{0\}} 
\tr \left(\Opw(a_{g,\lambda})\
\Opw \left[\FR (\varphi(g\,\cdot)) (\sqrt {|\lambda|} \, \cdot,\sqrt \lambda\, \cdot ,\lambda)\right]\right) 
|\lambda|^n d\lambda.
\end{equation}
For the value of the constant $c'_n$, see \cite{RF-monograph}.
The formula \eqref{EQ:quant2} now involves mainly  `Euclidean objects'.
 
\section{Some applications} 
 
We finally note several applications of the above theorems to questions of hypoellipticity
of (pseudo)\-differential operators on the Heisenberg group.
We say that a pseudo-differential operator $A$ is 
{\em Schwartz hypoelliptic} whenever 
$f\in \cS'(\bH_n), \ Af\in \cS(\bH_n)$ imply that $f\in \cS(\bH_n).$
Then, for example, as a simple consequence of our calculus we obtain that
the operator $\id-\cL$ is Schwartz hypoelliptic. In fact, criteria can be given in terms of the
$\lambda$-symbols:

\medskip

\begin{corollary} \label{thm_ellipticity_heis-cras}
Let $m\in \bR$ and $1\geq \rho\geq \delta\geq0$, $\rho\not=0$.
Let $\sigma=\{\sigma(g,\lambda)\}$ be in $S^m_{\rho,\delta}(\bH_{n})$ 
with $\sigma(g,\lambda) = \Opw \left( a_{g,\lambda} \right)$
as in Theorem \ref{thm_characterisation_Sm}.
Assume that  there are $R\in \bR$ and $C>0$ 
such that for any $(\xi,u)\in \bR^{2n}$ and $\lambda\not=0$ 
satisfying $|\lambda|(|\xi|^2+|u|^2)\geq R$ we have
\begin{equation}
\label{eq_a_elliptic_heis}
|a_{g,\lambda}(\xi,u)|\geq
C
\left(1+|\lambda|(1+|\xi|^2+|u|^2)\right)^{\frac{m} 2}.
\end{equation}
Then the operator $A$ in \eqref{eq_quantization} (or, alternatively, in \eqref{EQ:quant2})
has a left parametrix, i.e. there exists
$B\in \Psi^{-m}_{\rho,\delta}(\bH_{n})$ such that 
$BA-\id\in \Psi^{-\infty}$.
\end{corollary} 

\medskip
Corollary \ref{thm_ellipticity_heis-cras} has also a corresponding
`hypoellipticity version' which we omit here, but we give a few examples of both of them.
First, let $m,m_o\in 2\bN$ be two even integers such that $m\geq m_0$. 
Let $A$ be a differential operator given by either $X^m+iY^{m_o} +   T^{m_o/ 2}$
or $X^{m_o}+iY^m +   T^{m_o/ 2}$
on $\bH_1$. Then $A$ is Schwartz hypoelliptic 
and
satisfies the subelliptic estimates
$$
\forall s\in \bR
\quad\exists C>0\quad
\forall f\in \cS(\bH_1)\qquad
\|f\|_{L^p_{s+m_o}(\bH_{1})} 
\leq C\|A f\|_{L^p_{s}(\bH_{1})}.
$$
The above mentioned conclusion that $\id-\cL$ is Schwartz hypoelliptic can be also 
extended to variable coefficients using our calculus.
For example, if $f_1$ and $f_2$ are complex-valued smooth functions 
on $\bH_{n}$ such that 
$$
\inf_{x\in \bH_{n}, \lambda\geq \Lambda}
\frac{|f_1(x)+f_2(x)\lambda | }{1+\lambda} 
>0 \quad \mbox{for some} \ \Lambda\geq 0,
$$
and such that functions $X^{\alpha_1} f_1$, $X^{\alpha_2}f_2$
are bounded for every $\alpha_1,\alpha_2\in \bN_0^n$,
then
the differential operator $f_1(x)-f_2(x)\cL $ 
is Schwartz-hypoelliptic
and satisfies 
the following subelliptic estimates
$$
\forall s\in \bR
\quad\exists C>0\quad
\forall \varphi\in \cS(\bH_{n})\qquad
\|\varphi\|_{L^p_{s+2}(\bH_{n})} 
\leq C\|f_1 \varphi-f_2\cL \varphi\|_{L^p_{s}(\bH_{n})}.
$$

\section*{Acknowledgements}
The first author acknowledges the support of the London Mathematical Society via the Grace Chisholm Fellowship held at King's College London in 2011 as well as of the University of Padua.
The second author was supported in part by the EPSRC Leadership Fellowship EP/G007233/1
and both authors by EPSRC Grant EP/K039407/1.

\end{document}